\documentclass[a4paper,11pt]{article}
\usepackage[utf8]{inputenc}

\usepackage{amsmath}
\usepackage{amsfonts}
\usepackage{amssymb}

\usepackage[T1]{fontenc}
\usepackage{gensymb}
\usepackage{lmodern}
\usepackage{mdframed}

\usepackage[hyphens]{url}
\usepackage{hyperref}
\usepackage[hyphenbreaks]{breakurl}

\usepackage[left=3cm,right=3cm,top=3cm,bottom=3cm]{geometry}

\date{1 November 2018}
\author{
Maurice Chiodo\footnote{King's College, University of Cambridge. \texttt{mcc56@cam.ac.uk} .} \ \
  Dennis M\"uller\footnote{RWTH Aachen University \texttt{dennis.mueller3@rwth-aachen.de} .}   
}

\title{Mathematicians and Ethical Engagement}
\begin{document}

\maketitle

\let\thefootnote\relax\footnotetext{2020 \textit{AMS Classification:} 01A80. 00A30.}
\let\thefootnote\relax
\let\footnote\relax\footnotetext{Appeared in \textit{SIAM News 51, No. 9, 1 November 2018}. Available online at \url{https://sinews.siam.org/Details-Page/mathematicians-and-ethical-engagement}.}

\section*{}
In the past, some mathematical societies have discussed ethical policies and issues \cite{TheAmericanMathematicalSociety.2005, TheEuropeanMathematicalSociety.2012} and disseminated their own codes of conduct to address specific ethical concerns encountered by research mathematicians, such as those arising during publication. While ethical and behavioural issues specific to well-defined mathematical areas are of course still relevant, the last two decades have yielded many \textit{new} ethical concerns that now affect \textit{all} mathematicians in some way. Having taught these issues for more than two years at the University of Cambridge, we came to the realization that mathematicians can assume several different levels of ethical engagement \cite{Chiodo.2018}. Ethics in mathematics is not a binary process.

As the oldest consistently used scientific tool in Western thinking, mathematics carries perhaps the greatest scientific authority. It has become an extraordinarily powerful instrument ubiquitous to all of science and technology. How many hours of mathematical work underpin the technology behind smartphones, airplane flights, or models of global climate dynamics? But the applications—and therefore ethics—of mathematics go well beyond engineering. Modern mathematics is at the heart of economics and finance, and excessive trust in mathematical models contributed to the 2007 financial crisis. Even the most ardent purists in number theory or algebra can no longer claim to “just do the mathematics” and “leave the implications to ethicists”, as recent revelations about global mass surveillance have underscored their work‘s immediate social and political impact. It is now evident that one can wield practically \textit{all} branches of mathematics both for good and harm. Modern mathematics is a double-edged sword.

Just as physicists had to recognise the enormous ethical implications of their work after the atomic bombing of Hiroshima in August 1945, socially responsible mathematicians must also realise the existence of ethics in mathematical practise, which leads to issues far more complex and harder to characterize than publishing-related decisions. Plagiarism and the ethics of journal submission are real concerns, but hardly of the same order as these new ethical matters.

The inner workings of even areas of broad appeal—such as data science, machine learning, and optimisation—are often beyond the layman’s comprehension. Lawyers and judges struggle to understand policing and sentencing algorithms, politicians stretch to comprehend the full capabilities of state surveillance agencies, and electoral commissions barely grasp the algorithms and mathematical psychometrics behind Cambridge Analytica’s targeted advertising. Thus, only mathematicians can begin the process of unveiling the meaning, validity, applicability, and reliability of modern mathematics, paving the way for judges, politicians, and regulators to step in.

Even if we feel that mathematical research is beyond all ethical consideration, as academics we must ask ourselves: What do our students do after graduation? We train them in a wide range of mathematics, but do we teach them to be aware of possible ethical issues in its use? As a society, we have long agreed that the so-called Nuremberg defense—simply saying “I’m just doing my job” or “I was only following orders”—is not a valid excuse. Thus, it is imperative for us to teach ethics to our students and help them better contextualise their mathematical work. In April 2016, we began giving ethics seminars featuring guest speakers from industry, academia, and intelligence agencies to researchers and students in the Faculty of Mathematics at Cambridge. Shortly thereafter, we organized the first conference on “Ethics in Mathematics.”$^1$\let\thefootnote\relax\footnotetext{$^1$\url{http://www.ethics.maths.cam.ac.uk/EiM1/}} Through observation and case studies, we noticed that mathematicians can demonstrate what we term the “four levels of ethical engagement.” These levels form a recurring theme throughout our seminars.

The \textit{first} level is the fundamental understanding that the practice of mathematics is \textit{not} ethics free, and that ethical issues can surface in any mathematical work. One always performs mathematics in a social and political context, never in value-free isolation. Thus, all mathematicians must think about their individual responsibilities, as ethical issues may emerge at any time. This diligence can be as simple as considering environmental impact rather than merely optimising over time and money during a construction project. Mathematics can pose immediate or distant consequences that generally manifest as good, sometimes as not entirely good, and occasionally as downright bad. On this individualistic level, mathematicians modify and adapt their own ethical consciousness and actions, taking the important first step towards a more robust ethical awareness.

The \textit{second} of these four levels involves mathematicians \textit{speaking out} to other mathematicians, raising awareness of ethical issues among their peers. Individual mathematicians may recognize ethical issues in the mathematical work of others and try to inform them. They might precipitate unified action among their colleagues and locally bring about a collective ethical awareness and approach. Or they might write an article about ethics for their community, as we have done here.

The \textit{third} level is more complex. It teaches mathematicians to \textit{take a seat at the tables of power}. Mathematicians often need to learn the specific skills required to work with politicians, corporate management, and other non-scientists. These include engaging in policy discussions, establishing and rationalising their mathematical work’s objectives, and communicating potential limitations and possible drawbacks. Engineers and computer scientists are taught this at the undergraduate level, but mathematicians seldom receive such lessons explicitly. Many mathematicians in advancing industry careers unexpectedly find themselves in positions that require these abilities. Mathematics is becoming an increasingly powerful social tool, and seeing its creators hiding behind formulae and retrospectively apologising is not appropriate. If we want to take credit for our output’s positive impact, we should also be able to defend and properly contextualise our work and engage in apparently non-mathematical debates.

Our \textit{fourth} and final level is the responsibility of mathematicians to \textit{call out the bad mathematics of others} by proactively seeking out, learning about, and acting upon instances where mathematics has “gone wrong” — possibly in unrelated organisations. However, bad mathematics occurs in two distinct forms. First, it can refer to the practice of claiming results that are not mathematically true. The catastrophic misuse of statistics in the trial of Sally Clark \cite{Schneps.2013}, which the Royal Statistical Society reprimanded through the release of a statement \cite{TheRoyalStatisticalSociety.2001}, is one such example. Bad mathematics can also refer to trained mathematicians’ inappropriate use of mathematics by giving it excessive authority or directing it in ways that cause harm and exploit others. Members of any profession have the responsibility to hold their work—and the work of their colleagues—to high standards. Like statisticians, engineers, and doctors, mathematicians must adapt their own form of professional standards in academia, industry, and overall society. Some mathematicians are already questioning the validity and fairness of various decision-making algorithms or identifying the potential harms of artificial intelligence (AI), bringing such dangers into public consciousness and proposing workable solutions.

Practising ethics in mathematics is not binary, and mathematicians must consider various levels of engagement and ethical sensibility. Of course, our aforementioned four levels are an artificial and simplistic construct. One can refine them ad nauseum, but collectively they illustrate the depth and complexity of ethics in mathematics.

Not every mathematician will face problems pertaining to all levels, but everyone should remain aware of their social responsibilities, acknowledge the existence of ethical issues in the mathematical context, and appreciate their complexity. We teach students a broad spectrum of mathematics to prepare them for a wide variety of academic and professional eventualities. Why shouldn't we teach a broad spectrum of ethical situations in mathematics, which go beyond specialised courses such as ethics in AI? Lawyers, medics, biologists, engineers, physicists, and computer scientists learn subject-specific ethics because they will encounter these questions as professionals. Comprehending the seemingly-limitless uses of mathematics is difficult, and the ethical implications of modern mathematics depend on subtleties that only the mathematically-trained can understand. We are the only ones who can see behind the formulae. Thus, we should no longer leave these issues to professional ethicists and philosophers. No one else can address them, so we must.


\end{document}